\documentclass{amsart}
\usepackage{mathtools}
\usepackage{amsfonts,mathrsfs,dsfont}
\usepackage{graphicx,xcolor,hyperref} 

\newtheorem{theorem}{Theorem}

\theoremstyle{definition}

\theoremstyle{remark}
\newtheorem*{remark}{Remark}

\def \R {\mathbb{R}}

\def \T {\mathcal{T}}

\DeclareMathOperator\supp{supp}
\DeclareMathOperator{\sgn}{sgn}
\newcommand{\overbar}[1]{\mkern 1.5mu\overline{\mkern-1.5mu#1\mkern-1.5mu}\mkern 1.5mu}

\DeclarePairedDelimiter \mscal{\langle}{\rangle_M}

\numberwithin{equation}{section}

\title{Recovery of a Null Form in the Wave Equation from Scattering Data}
\author{Joel Nathe}
\date{}

\address{\parbox{\linewidth}{
  Department of Mathematics, Purdue University\\
150 North University Street, West Lafayette IN 47907, USA}}
\email{jnathe@purdue.edu}

\begin{document}
\begin{abstract}
    We use highly oscillatory geometric optics solutions to solve the inverse problem for the system
    $$
    \square \begin{bmatrix} u^{(1)}\\ u^{(2)}\\ \vdots\\ u^{(n)}\end{bmatrix} = \sum_{k,l=1}^n \left(Q_0(u^{(k)},u^{(l)}) \begin{bmatrix} q_{1kl}\\q_{2kl}\\\vdots\\q_{nkl} \end{bmatrix} + \sum_{i>j=1}^d  Q_{ij}(u^{(k)},u^{(l)})\begin{bmatrix} p_{1ijkl} \\ p_{2ijkl} \\ \vdots \\ p_{nijkl} \end{bmatrix}\right), 
    $$
    where $Q_0$ and $Q_{ij}$ are the symmetric and anti-symmetric bilinear null forms. We present solutions in both the linear and weakly nonlinear regimes. In the linear regime, we show that the coefficients of order $h^3$ determine an injective light-ray transform of a vector field which depends on the coefficients $q_{rkl}, p_{rijkl}$. In the weakly nonlinear regime, we see that the coefficients of order $h$ determine the non-abelian light ray transform for matrices associated with the coefficients. While we do not have an injectivity result for this case, we do have one if we assume the coefficients do not depend on the time variable $x_0$, as our coefficients instead determine an injective non-abelian X-ray transform.
\end{abstract}
\maketitle
\section{Introduction}
We consider the inverse problem for a system of wave equations with a nonlinearity that satisfies what is called a null condition. Specifically, if we write $x \in \R^{1+d}$ ($d \geq 2$) as $(x_0,x')$ for $x' \in \R^d$, and define the null forms
\begin{align*}
Q_0(u^{(k)},u^{(l)}) &\coloneqq -\partial_{x_0} u^{(k)} \partial_{x_0} u^{(l)} + \sum_{j=1}^d \partial_{x_j} u^{(k)} \partial_{x_j} u^{(l)} \\
Q_{ij}(u^{(k)},u^{(l)}) &\coloneqq \partial_{x_i} u^{(k)}\partial_{x_j} u^{(l)} - \partial_{x_j} u^{(k)} \partial_{x_i} u^{(l)},
\end{align*}
our goal is to solve the inverse problem for the system
\begin{equation}\label{eq.system}
\square \begin{bmatrix} u^{(1)}\\ u^{(2)}\\ \vdots\\ u^{(n)}\end{bmatrix} = \sum_{k,l=1}^n \left(Q_0(u^{(k)},u^{(l)}) \begin{bmatrix} q_{1kl}\\q_{2kl}\\\vdots\\q_{nkl} \end{bmatrix} + \sum_{\substack{i,j=1\\i>j}}^d  Q_{ij}(u^{(k)},u^{(l)})\begin{bmatrix} p_{1ijkl} \\ p_{2ijkl} \\ \vdots \\ p_{nijkl} \end{bmatrix}\right)
\end{equation}
where $\square$ is the d'Alembertian $-\partial_t^2 + \sum_{j=1}^d \partial_{x_j}^2$ and we solve for $q_{rkl}(x,u)\in C_0^\infty(\R^{1+d} \times \R^n)$ and $p_{rijkl}(x,u) \in C_0^\infty(\R^{1+d}\times \R^n)$. We also assume $p_{rijlk} = -p_{rijkl}$ and $q_{rkl} = q_{rlk}$, since $Q_0$ being symmetric and $Q_{ij}$ anti-symmetric implies that \ref{eq.system} depends only on $q_{rkl}+q_{rlk}$ and $p_{rijkl}-p_{rijlk}$. And for notation purposes, we set $p_{rijkl} = -p_{rjikl}$, let
\begin{equation*}
    \begin{gathered}
    \nabla_{x} u \coloneqq  (\partial_{x_0}u,\partial_{x_1} u,...,\partial_{x_d} u), \\
    \nabla_{x'}u \coloneqq (\partial_{x_1} u,...,\partial_{x_d} u),
    \end{gathered}
\end{equation*} 
and define the operator $Q_Tu$ to be the right-hand side of \eqref{eq.system}.\par
\par
The short-time existence and uniqueness of the Cauchy problem for the system has been studied by Sogge \cite{Sogge}. However, we will use a method involving highly oscillatory solutions which does not need this result, as in solving the inverse problem we will not need to consider arbitrary initial data. This method has been used by the author and S\`a Barreto in \cite{NSB} to treat the scalar ($n=1$) case of this problem. We will extend these results to the general case. 

The null forms we consider arise when considering the regularity needed for a wave equation initial value problem $\square u = F(u,\nabla u)$ to be well-posed. If the nonlinearity $F$ is of the form of $Q_T$, then the regularity of the initial data needed is less than the general case \cite{KM1,Sogge}. The imposition of a null condition can also lead to global existence results, as in \cite{Klainerman,Sideris}.
\par
These regularity and existence improvements are relevant to the various physical equations in which null forms occur. The most notable of these is the wave maps equation, which when put into the local coordinates of a Riemannian manifold (with the convention that we sum over repeated indices) is
\[
\square u^r = \Gamma_{jk}^r(u) \partial^\alpha u^j \partial_\alpha u^k,
\]
where $\partial^{x_i}= \partial_{x_i}$ for $i > 0$ and $\partial^{x_0} = -\partial_{x_0}$ \cite{Tataru}. This is \eqref{eq.system} when all $p_{rijkl}$ terms are zero---that is, with only the $Q_0$ forms. 
Null forms also appear in the Einstein field equations \cite{LR} and Yang-Mills equations \cite{Christodoulou}.
\par
To solve our inverse problem, we will use highly oscillatory solutions, known as nonlinear geometric optics solutions. The techniques of nonlinear geometric optics have a long history in physics \cite{Bloembergen, BW, Boyd}, but while some work was done by Lax in the 1950s \cite{Lax1,Lax2} the mathematics were not rigorously set out until the late 1980s and early 1990s, with work by Gu\`es \cite{Gues}, Majda \cite{Majda}, and the collaboration of Joly, M\'etivier, and Rauch \cite{JMR1,JMR2,JMR3}. The use of these techniques for inverse problems involving the recovery of wave equation nonlinearities began with \cite{SBP1}, and has continued in \cite{ES1,ES2,NSB,SBP2}. 
\par
The interest in these inverse problems for nonlinear hyperbolic equations has grown following the publication of \cite{KLU} by Kurylev, Lassas, and Uhlmann. The techniques used in \cite{KLU} as well as the subsequent works \cite{FLO, HU, KLOU,LUW, UZ1,UZ2} involve the propagation of singularities, and generally aim to recover information about a manifold. While our methods may achieve that result, such as in the case of the wave maps equation, they are instead focused on recovering the nonlinearity.
\par
We give two possible ways to solve the inverse problem, one using the linear regime, and the other using the weakly nonlinear regime (see \cite{MetivierBook}, chapter 5). The existence of oscillatory solutions for nonlinearities satisfying such null conditions has been studied \cite{Sogge}, but the above inverse problem is new, apart from the $n=1$ case previously treated \cite{NSB}.
\par
In both regimes, we analyze the oscillatory solution at a given time $T$.
For the linear regime, we look at the coefficient of order $h^3$ to obtain a light ray transform which we know to be injective. For the weakly nonlinear regime, we instead use the coefficient of order $h$. However, to obtain results from observations of data of this order, we must sacrifice some of our conditions. For this transform, we only know the X-ray version (where our unknown coefficients are independent of time) to be injective.
\par
If we have $V = (1,\theta) \in \R^{1+d}, \theta \in \mathbb{S}^{d-1}$, we write
\begin{align*}
    \widetilde{V} &= (-1,\theta)\\
    \langle x,V\rangle_M &= -x_0 + \sum_{j=1}^d x_j \theta_j\\
    \T_{V} u &= \partial_{x_0} u + \theta \cdot \nabla_{x'} u
\end{align*}
where $\langle x,V\rangle_M$ is the Minkowski scalar product, and $\T_{V}$ is a transport operator.
If $\varphi \in C^\infty(\R)$, we denote 
\begin{equation}\label{eq.background_definition}
\varphi_V(x) = \varphi(\langle x,V\rangle_M),\quad \varphi'_V(x) = \varphi'(\langle x,V\rangle_M),\quad \varphi''_V (x) = \varphi''(\mscal{x,V})
\end{equation}
and notice that
\begin{equation}\label{eq.backgroundsol}
    \square \varphi_V = \T_V \varphi_V = Q_0(\varphi_V,\varphi_V) = Q_{ij}(\varphi_V,\varphi_V) = 0.
\end{equation}
\par
%
In the linear regime, our oscillatory solutions will be of the form
\begin{equation}
u(x) = \begin{bmatrix} u^{(1)}(x)\\u^{(2)}(x)\\\vdots\\u^{(n)}(x)\\ \end{bmatrix} = C+h\varphi_V(x)\vec{e}_l + h^2\sum_{\pm m=0}^\infty e^{\frac{im}{h}\mscal{x,W}} \begin{bmatrix}  A_m^{(1)}(h,x)\\  A_m^{(2)}(h,x)\\ \vdots\\  A_m^{(n)}(h,x)\end{bmatrix}
\end{equation}
where $\vec{e_l} \in \R^n$ is the $l$th standard basis vector and $C \in \R^n$ is a constant.\par

Due to \eqref{eq.backgroundsol}, we can think of $C+\varphi_V(x) \vec{e}_l$ as a background solution, which we perturb with oscillatory expansions. That this view holds for all choices of $V$ and $l$ gives us the flexibility to solve for all coefficients of our inverse problem.
\par
Specifically, we perturb $C+h\varphi_V(x) \vec{e}_l$ with $h^2(\frac{A\pm iB}{2}) e^{\mp\frac{i}{h}\mscal{x,W}} \vec{e_k}$. While the nonlinearity vanishes for each of these terms individually, their interaction with each other results a term of order $h^2$, which must be matched in the solution by an oscillatory term of order $h^3$. It is this term that we measure to find the light ray transform of our coefficients, and thus solve the inverse problem. If $l=k$, this term depends solely on the $q_{rkl}$ coefficients. If instead $l \neq k$, it also depends on the $p_{rijkl}$ coefficients.
\par
Of course, since our equation is nonlinear, each new term ``produces" more terms itself. We show that as $h\to0$, these later terms decay faster than the terms of order $h^3$ we are interested in. We also use this decay to show that such an oscillatory solution exists.
\begin{theorem} \label{thm.main}
    Assume all $q_{rkl},p_{rijkl}$ are supported on $\{|x| \leq R \} \times \R^n$. Let $\varphi,\chi \in C_0^\infty(\R)$ be real-valued, $V = (1,\theta), W = (1,\omega)$ with $\theta,\omega \in \mathbb{S}^{d-1}$, and $\varphi_V,\chi_W$ defined as in \eqref{eq.background_definition}. Now choose $k,l\in \mathbb{Z}$ with $1 \leq k,l \leq n$, $k \geq l$. For $A,B \in \R$, $C\in \R^n$, consider 
    \begin{equation}\label{eq.incoming_wave}
    u_{inc}(h,V,W,C,x) = C+h\varphi_V(x) \vec{e_l}
     +h^2\chi_W(x)\left(A \cos\left(\frac{\mscal{x,W}}{h} \right) + B \sin\left(\frac{\mscal{x,W}}{h} \right)\right) \vec{e_k}
    \end{equation}
    where $\vec{e_l},\vec{e_k}$ are the $l$th and $k$th standard basis vectors, respectively, in the solution space $\R^n$. This satisfies \eqref{eq.system} when $x \notin \supp \chi_W$. 
    \par
    Then for any integer $N \geq 2$ and any interval $[T_0,T] \times \R^d$, there exists a unique solution $u$ of \eqref{eq.system} with
    \[
    u(h,V,W,C,x)= u_{inc}(h,V,W,C,x)
    \]
    for all $x$ with $x_0 \ll 0$ outside the support of $q_{rkl},p_{rijkl}$, such that $u$ has an expansion of the form 
    \begin{multline}
    u(h,V,W,C,x) = C + h\varphi_V(x)\vec{e_l} + h^2\sum_{\pm m=0}^N e^{\frac{im}{h}\mscal{x,W}} \begin{bmatrix}  \sum_{p=0}^{N-2} h^pA_{m,p}^{(1)}(x)\\  \sum_{p=0}^{N-2} h^pA_{m,p}^{(2)}(x)\\ ...\\ \sum_{p=0}^{N-2} h^p A_{m,p}^{(n)}(x)\end{bmatrix}\\ + h^{N+1} \mathcal{E}_N(h,V,W,C,x).
    \end{multline}
    \par
    Additionally, the coefficients satisfy $A_{m,p}^{(r)} \in C^\infty(\R^{1+d}), \overbar{A_{m,p}^{(r)}} = A_{-m,p}^{(r)}$, so $u$ is smooth and real-valued. And for $x_0 \ll 0$, the coefficients satisfy the initial conditions
    \begin{equation}\label{eq.initial_conditions_abel}
        A_{m,p}^{(r)} =  \begin{cases}
            \frac{1}{2} \chi_W(x)(A-iB\sgn{m}) & r=k,|m|=1,\text{ and } p=0\\
            0 & \text{otherwise}.
        \end{cases}
    \end{equation}
    For this solution, the above describes $A_{m,0}^{(r)}$ for all $x$. In addition, for any $\tau > 0$, $k \in \mathbb{N}$, the error term $\mathcal{E}_N$ satisfies
    \begin{equation}\label{eq.exact_error_bounds_abel}
    \begin{gathered}
    \mathcal{E}_N,\nabla_x \mathcal{E}_N \in C^0 ([T_0,T], H^m(\R^d_{x'})),\\
    \sup_{x_0 \in [-\tau,\tau]} || \nabla_{x'}^k \mathcal{E}_N(x_0,\cdot)||_{H^1(\R^d)} + \sup_{x_0 \in [-\tau,\tau]} || \nabla_{x'}^k \partial_{x_0} \mathcal{E}_N(x_0,\cdot)||_{L^2(\R^d)} \leq C_k(N,\tau)h^{-k-1}
    \end{gathered}
    \end{equation}
    for sufficiently small $h$.
    \par
    If $k = l$, the coefficients $A_{1,1}^{(r)}$ satisfy the transport equation  
    \begin{equation}\label{eq.transport_sym_abel}
    \T_W A_{1,1}^{(r)} = q_{rll}(x,C)F(V,W,x) A_{1,0}^{(l)}.
    \end{equation}
    where
    \[
    F(V,W,x) = \varphi'_V(x) \langle V,W\rangle_M.
    \]
    In this case, for $T$ large enough, $u(h,V,W,x)|_{x_0 = T}$ completely determines $q_{rll}$ for all $r$.
    \par
    If instead $k \neq l$, the coefficients $A_{1,1}^{(r)}$ satisfy 
    \begin{align}\label{eq.transport_anti_abel}
        \T_W A_{1,1}^{(r)} &= \varphi_V'(x) \left( q_{rkl}(x,C)\mscal{V,W} +\sum_{\substack{i,j = 1\\i > j}} p_{rijkl}(x,C) (V_iW_j -V_j W_i)\right) A_{1,0}^{(k)}\\ 
        &=\langle G(k,x,V,C),W\rangle_M A_{1,0}^{(k)}
    \end{align}
    where
    \[
    G(r,x,V,C) \coloneqq \varphi'_V(x)\begin{bmatrix} q_{rkl}(x,C) &  0 & 0 & ... & 0\\
    0 & q_{rkl}(x,C) & p_{r21kl}(x,C) &... & p_{rn1kl}(x,C)\\ 0 & -p_{r21kl}(x,C) & q_{rkl}(x,C)  & ... &  p_{rn2kl}(x,C)\\
    \vdots & \vdots& \vdots&  \ddots & \vdots\\
    0 & -p_{rn1kl}(x,C) & -p_{rn2kl}(x,C)  & ... & q_{rkl}(x,C)
    \end{bmatrix}V.
    \]
    For $T$ large enough, $u(h,V,W,C,x)|_{x_0=T}$ completely determines $p_{rijkl}$ for all $r,i,j$.
\end{theorem}
\begin{remark}
    While our result requires that $q_{rkl}, p_{rijkl}$ are bounded in the time variable $x_0$, in practice one can usually solve for a nonlinearity supported in a cylinder $\{|x'| < R\} \times \R^n$ by applying a cutoff function.
\end{remark}
We can also work in the weakly nonlinear regime. Here, our background solution $\varphi_V \vec{e_l}$ will have order $h^0$, and our oscillatory expansions will have order $h^1$. The terms produced by the interaction between these are the same order as the first transport equation, so we obtain matrix valued transport equations.
\begin{theorem}\label{thm.non_abelian}
    Under the same assumptions as Theorem \ref{thm.main}, except that $q_{rkl},p_{rijkl}$ are supported on the cylinder $\{|x'| < R\} \times \R^n$, consider
    \begin{equation}\label{eq.incoming_wave_nonabel}
        v_{inc}(h,V,W,x) = \varphi_V(x) \vec{e_l} + h \chi_W(x) \left(A \cos\left(\frac{\mscal{x,W}}{h} \right) + B \sin\left(\frac{\mscal{x,W}}{h} \right)\right) \vec{e_k}.
    \end{equation}
    As in Theorem \ref{thm.main}, this satisfies \eqref{eq.system} when $x \notin \supp \chi_W$. 
    \par
    Then for any positive integer $N$ and any interval $[T_0,T] \times \R^d$, there exists a unique solution $v$ of \eqref{eq.system} such that $v(h,V,W,x) = v_{inc}(h,V,W,x)$ for all $x \ll 0$ with $x$ outside the support of $q_{rkl},p_{rijkl}$,
    and $v$ has an expansion of the form
    \[
    v(h,V,W,x) = \varphi_V(x) \vec{e_l} + h\sum_{\pm m = 0}^{N} e^{\frac{im}{h}\mscal{x,W}} \chi_W(x)\begin{bmatrix} \sum_{p=0}^{N-1} h^p B_{m,p}^{(1)}\\   \sum_{p=0}^{N-1} h^p B_{m,p}^{(2)}\\
    ...\\
     \sum_{p=0}^{N-1} h^p B_{m,p}^{(n)} \end{bmatrix} + h^{N+1} \tilde{\mathcal{E}}_N(h,V,W,x).
    \]
    As in Theorem \ref{thm.main}, the coefficients satisfy $B_{m,p}^{(r)} \in C^\infty(\R^{1+d}), \overbar{B_{m,p}^{(r)}} = B_{-m,p}^{(r)}$, so $v$ is smooth and real-valued. And for $x_0 \ll 0$, the coefficients satisfy the initial conditions
    \begin{equation}\label{eq.initial_conditions_nonabel}
        B_{m,p}^{(r)} =  \begin{cases}
            \frac{1}{2} \chi_W(x)(A-iB\sgn{m}) & r=k,|m|=1,\text{ and } p=0\\
            0 & \text{otherwise.}
        \end{cases}
    \end{equation}
    For this solution, the error term $\tilde{\mathcal{E}}_N$ satisfies the bounds given in \eqref{eq.exact_error_bounds_abel} for $\mathcal{E}_N$.
    \par
    In addition, the coefficient vector $B_{1,0}$ satisfies the equation
    \begin{equation}\label{eq.transport_nonabel}
    \T_W B_{1,0} = \Phi B_{1,0} + \sum_{i=1}^d \omega_i \Psi_i B_{1,0}
    \end{equation}
    where
    \begin{align}
        \Phi &= \begin{bmatrix} q_{11l} & q_{12l}& ... & q_{1nl}\\ q_{21l}& q_{22l} & ...&q_{2nl}\\ \vdots  &\vdots& \ddots & \vdots\\ q_{n1l} & q_{n2l} & ... & q_{nnl} \end{bmatrix}    \\
        \Psi_i &= -\theta_i \Phi + \sum_{j=1}^d \theta_j
        \begin{bmatrix}
            p_{1ji1l} & p_{1ji2l} & ... & p_{1jinl}\\ p_{2ji1l} & p_{2ji2l} & ... & p_{2jinl}\\
            \vdots & \vdots & \ddots & \vdots\\
            p_{nji1l} & p_{nji2l} & ... &p_{njinl}.
        \end{bmatrix}
    \end{align}
    If $q_{rkl},p_{rijkl}$ are independent of $x_0$, then $v(h,V,W,x)|_{x_0 = T}$ for sufficiently large $T$ completely determines all $q_{rkl},p_{rijkl}$.
\end{theorem}
\section{The Inverse Problem} We prove the end of Theorem \ref{thm.main}, that $u(h,V,W,x)$ evaluated at $x_0 = T$ large enough completely determines our coefficients. We show that due to the transport equations \eqref{eq.transport_sym_abel} and \eqref{eq.transport_anti_abel}, the evaluation determines an injective light ray transform of an object dependent on the coefficients.
\par
We also prove the final claim of Theorem \ref{thm.non_abelian}, that if $q_{rkl}, p_{rijkl}$ are independent of $x_0$, then due to \eqref{eq.transport_nonabel}, the evaluation of $v(h,V,W,x)$ at $x_0 = T$ large enough determines the coefficients. 
\par
\subsection{The Coefficients $q_{rll}$}
We start with the case $k=l$, where we focus on the $q_{rll}$ coefficients. If we let $x' = y+x_0 \omega$, $x_0 = s$, in the new coordinates \eqref{eq.transport_sym_abel} becomes
\[
\partial_s A_{1,1}^{(r)}(s,y) = q_{rll}(s,y,C) F(V,W,s,y) A_{1,0}^{(l)}(s,y).
\]
\par
The initial condition is $A_{1,1}^{(r)} = 0$ for $s \ll 0$. Solving this equation gives
\[
A_{1,1}^{(r)}(s,y) = \int_{-\infty}^s A_{1,0}^{(l)}(t,y) q_{rll}(t,y,C) F(V,W,t,y) dt.
\]
Changing back to $x_0,x'$,
\[
A_{1,1}^{(r)}(x_0,x') = \int_{-\infty}^{x_0} A_{1,0}^{(l)}(t,x'+(t-x_0)\omega)  q_{rll}(t,x'+(t-x_0)\omega,C) F(V,W,t,x'+(t-x_0)\omega) dt
\]
We can shift the time variable so that $q_{rll}$ has support contained in $\{x_0 < 0\}$ (i.e. we take our measurement at $x_0 = T = 0$). Then the above equation becomes
\begin{align*}
A_{1,1}^{(r)}(0,x') &= \int_{\R}A_{1,0}^{(l)}(t,x'+t\omega) q_{rll}(t,x'+t\omega,C) F(V,W,t,x'+t\omega) dt\\
&= \int_{\R}  \langle  \mathscr{F}(q_{rll},A_{1,0}^{(l)} ,V,C,t,x'+t\omega), W\rangle dt
\end{align*}
where $\langle\cdot,\cdot\rangle$ is the usual $\R^{d+1}$ inner product, and 
\[
\mathscr{F}(q_{rll},A_{1,0}^{(l)},V,C,t,x'+t\omega) = A_{1,0}^{(l)}(t,x'+t\omega)q_{rll}(t,x'+t\omega,C) \varphi_V'(t,x'+t\omega) \widetilde{V}.
\]
The above integral is the future light-ray transform of $\mathscr{F}(q_{rll},A_{1,0}^{(l)},V,C,\cdot,\cdot)$, defined as
\[
L_1(\mathscr{F})(x',W) = \int_{\R} \langle  \mathscr{F}(t,x'+t\omega), W\rangle dt.
\]
\par
So we expand on the conclusion of Theorem \ref{thm.main} to see that the solution is of the form 
\begin{multline}\label{eq.sym_order_3_l}
u^{(l)}(h,V,0,x') = h \varphi_V(x) + h^2 \chi_W(x) \left( A  \cos\left(\frac{i}{h} \langle x,W\rangle_M\right) + B\sin\left(\frac{i}{h} \langle x,W\rangle_M\right)\right)\\
+ h^3 \Bigg(2\Re {L_1(\mathscr{F}(q_{lll},A_{1,0}^{(l)},V,C,x',W))}\cos\left(\frac{\mscal{x,W}}{h} \right) \\+ 2\Im {L_1(\mathscr{F}(q_{lll},A_{1,0}^{(l)},V,C,x',W))} \sin\left(\frac{\mscal{x,W}}{h} \right)\Bigg)+ O(h^4)
\end{multline}
and
\begin{multline}\label{eq.sym_order_3_r}
   u^{(r)}(h,V,0,x') =  h^3 \Bigg(2\Re {L_1(\mathscr{F}(q_{rll},A_{1,0}^{(l)},V,C,x',W))}\cos\left(\frac{\mscal{x,W}}{h}\right) \\+ 2\Im L_1(\mathscr{F}(q_{rll},A_{1,0}^{(l)},V,C,x',W)) \sin\left(\frac{\mscal{x,W}}{h} \right)\Bigg) + O(h^4).
\end{multline}
for $r \neq l$.
Since $L_1$ is linear, $A_{1,0}^{(l)}$ is known and $\chi,\varphi,C$ are arbitrary, showing that \eqref{eq.sym_order_3_l} and \eqref{eq.sym_order_3_r} determine $q_{rll}$ for all $r$ is equivalent to showing that if $L_1(\mathscr{F}) = 0$, $\mathscr{F} = 0$.
\par
It is known that the light ray transform $L_1$ is injective modulo potential fields (see Section III.4 of \cite{SU}, as well as \cite{RabienHaratbar,Stefanov}), so it suffices to show that $d\mathscr{F} = 0$ implies $\mathscr{F} = 0$. This has already been done in the previous paper solving the $n=1$ case (see Proposition 1.3 of \cite{NSB}). Thus, $u(h,V,W,x)|_{x_0=T}$ completely determines $\mathscr{F}$. 
\par
While the light ray transform is invertible, the inverse is unstable. However, if all $q_{rll}$ are independent of $x_0$, with some adjustment to $\varphi_V$ and $\chi_W$ the light-ray transform becomes a stably invertible X-ray transform \cite{SU}.
\subsection{The Coefficients $q_{rkl}$, $p_{rijkl}$}
We now consider the case $k \neq l$. If we let $x' = y+x_0\omega,x_0=s$, in the new coordinates \eqref{eq.transport_anti_abel} becomes 
\[
\partial_s A_{1,1}^{(r)}(s,y) = \langle G(s,y,V,C),W\rangle_M A_{1,0}^{(l)}.
\]
With initial condition $A_{1,1}^{(r)}= 0$ for $s = 0$, this has solution
\begin{align*}
A_{1,1}^{(r)}(s,y) &= \int_{-\infty}^s A_{1,0}^{(k)}(t,y) \langle G(k,t,y,V,C),W\rangle_M dt.\\
A_{1,1}^{(r)}(x_0,x')&= \int_{-\infty}^{x_0} A_{1,0}^{(k)}(t,x'+(t-x_0)\omega) \langle G(k,t,x'+(t-x_0)\omega,V,C),W\rangle_M dt.
\end{align*}
If we shift the time variable appropriately, we have
\begin{equation}\label{eq.anti_light_ray}
A_{1,1}^{(r)}(0,x') = \int_{\R} \langle \mathscr{G}(k,A_{1,0}^{(k)},V,C,t,x'+t\omega),W\rangle_M dt = L_1(\mathscr{G}(k,A_{1,0}^{(k)},\tilde{V},C,\cdot,\cdot))(x',W)
\end{equation}
where $\mathscr{G}(k,A_{1,0}^{(k)},V,t,x'+t\omega) = A_{1,0}^{(k)}(t,x'+t\omega) G(k,t,x'+t\omega,V,C)$.
\par
Following the argument from the $k=l$ case, we must show that $d\mathscr{G} =0$ implies that $q_{rkl},p_{rijkl} = 0$. We see that
\begin{multline}
    d\mathscr{G} = \sum_{i=1}^d \Bigg(\Big((\partial_{0}-\partial_{i})(A_{1,0}^{(k)} q_{rkl})\varphi_V' - 2\varphi_V'' \theta_i A_{1,0}^{(k)} q_{rkl}\\
    + \sum_{j=1}^n \theta_j( \partial_{x_0}(A_{1,0}^{(k)} p_{rjikl})\varphi_V' - \varphi_V'' A_{1,0}^{(k)} p_{rjikl}) dx^0 \wedge dx^i\Big)\\
    + \sum_{\substack{m = 1\\ m< i}}^d \Big((\theta_i \partial_m(A_{1,0}^{(k)}q_{rkl}) - \theta_m \partial_i(A_{1,0}^{(k)}q_{rkl})) \varphi_V'\\
    + \sum_{j=1}^d \theta_j \Big( (\partial_m (A_{1,0}^{(k)} p_{rjikl}) - \partial_i (A_{1,0}^{(k)}p_{rjmkl})) \varphi_V' + (\theta_m p_{rjikl} - \theta_i p_{rjmkl})A_{1,0}^{(k)} \varphi_V''\Big) dx^m \wedge dx^i \Bigg)
\end{multline}
Since our choice of $\varphi$ was arbitrary, looking at $dx^0 \wedge dx^i$ implies that we must have $2\theta_i A_{1,0} q_{rkl} + \sum_{j=1}^d \theta_j A_{1,0}^{(k)} p_{rijkl} = 0$. Since $p_{riikl} = 0$ and $\theta, \chi$ were arbitrary, we must have $q_{rkl} = p_{rijkl} =0$ for all $i,j$.
\par
Since $k,l$ (with $k \geq l$) were arbitrary, the results of Theorem \ref{thm.main} show that $u|_{x_0=T}$ determines all $q_{rkl},p_{rijkl}$.

\subsection{The weakly nonlinear regime and the non-Abelian X-Ray transform}
For the inverse problem of Theorem \ref{thm.non_abelian}, if we set $x' = y+x_0 \omega$,  $x_0 = s$, \eqref{eq.transport_nonabel} is now
\begin{equation}\label{eq.connection_higgs}
\partial_s B_{1,0} = \Phi B_{1,0} + \sum_{i=1}^d \omega_i \Psi_i B_{1,0}.
\end{equation}
With a proper choice of $\varphi$, we can assume $\Phi,\Psi_i$ are independent of $x_0$. If $\Phi \equiv 0$, the above is the equation for the parallel transport of a connection. The $\Phi$ term added is known as a Higgs field \cite{PSU}. 
\par
If we replace $B_{1,0}$ with a matrix $U$, we call the solution to \eqref{eq.connection_higgs} with $\lim_{s \to -\infty} U = I$ the scattering data of $\Phi,\Psi$. Then if $B_{1,0}$ solves \eqref{eq.connection_higgs} with $\lim_{s\to-\infty} B_{1,0} = b$ for some initial data $b$, we have that $B_{1,0}(y,s) = U(y,s)b$. Since $\Phi,\Psi$ are compactly supported in $x'$, by taking $s = T$ large enough we obtain $C_{\Phi,\Psi} b$, where $C_{\Phi,\Psi} = \lim_{s\to\infty} U$ is the non-Abelian X-Ray transform of $\Phi,\Psi$.
\par
The non-Abelian X-Ray transform is known to be invertible up to gauge (see Chapter 14 of \cite{PSU}). This means that if $C_{\Phi,\Psi} = C_{\tilde{\Phi},\tilde{\Psi}}$ for all $(y,\omega)$, we must have
\begin{align*}
\tilde{\Phi} &= g \Phi g^{-1}\\
\tilde{\Psi}_i &= g \Psi_i g^{-1} + g^{-1} \partial_i g
\end{align*}
for some matrix $g(x')$ with $\lim_{|x'| \to \infty} g(x') = I$.
\par
However, we know that when $\omega=\theta$, we should obtain $\Phi+\sum_{i=1}^d \omega_i \Psi_i = 0$. This can only be achieved when $g \equiv I$ on the support of $\Phi,\Psi$, so our data does in fact completely determine $\Phi,\Psi$. 
\par
Then since $\theta$ and $l$ were arbitrary, we see that $v(h,V,W,x)|_{x_0 = T}$ completely determines all $q_{rkl}, p_{rijkl}$.
\begin{remark}
    If $q_{rkl}, p_{rijkl}$ depend on $x_0$, then the problem is that of the non-Abelian light ray transform. In \cite{Zyskin}, Zyskin finds invertibility for a related problem. Zyskin also states that this transform is not invertible, but does not give proof. We are not aware of any results related to the injectivity of this transform.
    \par
    If instead we know a priori that $q_{rkl},p_{rijkl} = 0$ when $k \neq l$ (so we have only $Q_0(u^{(l)},u^{(l)})$ terms), \ref{eq.transport_nonabel} reduces to a collection of abelian transport equations. We can show that these determine all $q_{rll}$ in a manner similar to that of the linear regime; this is what occurs in \cite{NSB} (the $n=1$ case).
\end{remark}
\section{Approximate Solutions}
We now construct approximate solutions of the form described in Theorem \ref{thm.main}, which we will later prove converge to the a solution of \eqref{eq.system}.
\begin{theorem}
    \label{thm.approx}
    Let $p_{rijkl},q_{rkl},\varphi,\chi,V,W,C,k,l$ be as in Theorems \ref{thm.main},\ref{thm.non_abelian} and let $u_{inc}$ be defined as in \eqref{eq.incoming_wave}, $v_{inc}$ defined as in \eqref{eq.incoming_wave_nonabel}. Then for any positive integer $N$, there exists approximate solutions $u_N,v_N$ to \eqref{eq.system}, in the sense that
    \begin{equation}\label{eq.approx_eval}
    \begin{gathered}
        \square u_N - Q_T u_N = h^{N+1}E_N(h,V,W,x)\\
        \square v_N - Q_T v_N = h^{N+1}\tilde{E}_N(h,V,W,x)
        \\
        u_N = u_{inc}, v_N = v_{inc} \text{ for $x_0 \ll 0$ outside the support of $p_{rijkl},q_{rkl},$}\\
    \end{gathered}
    \end{equation}
    and $E_N,\tilde{E}_N$ are smooth and satisfy \ref{eq.exact_error_bounds_abel}. Also, for any $\tau > 0$, $k \geq 1$, $u_N$ and $v_N$ satisfy
    \begin{equation}\label{eq.approx_error_est}
    \begin{gathered}
    u_N,\nabla_x u_N \in C^0([T_0,T],W^{k,\infty}(\R^d_{x'})),
    \sup_{x_0 \in [-\tau,\tau]} ||u_N(x_0,\cdot)||_{W^{1,\infty}(\R^d)} \leq C(N,\tau)\\
    \sup_{x_0 \in [-\tau,\tau]} ||\nabla_{x'}^k u_N(x_0,\cdot)||_{L^\infty(\R^d)}\leq C(N,\tau) h^{1-k},\sup_{x_0\in[-\tau,\tau]} ||\nabla_{x'}^k \nabla_x u_N(x_0,\cdot)||_{L^\infty(\R^d)} \leq C(N,\tau)h^{-k}.
    \end{gathered}
    \end{equation}
    Moreover, $u_N$ has an expansion of the form
    \[
    u_N(h,V,W,x) = C+h\varphi_V(x)\vec{e_l} + h^2\sum_{\pm m=0}^N e^{\frac{im}{h}\mscal{x,W}} \begin{bmatrix}  \sum_{p=0}^{N-2} h^pA_{m,p}^{(1)}(x)\\ \sum_{p=0}^{N-2} h^pA_{m,p}^{(2)}(x)\\ ...\\ \sum_{p=0}^{N-2} h^p A_{m,p}^{(n)}(x)\end{bmatrix}
    \]
    and is unique up to order $N+1$---that is, if $w_N$ satisfies \eqref{eq.initial_conditions_abel}, then
    \[
    \sup_{x_0 \in [-\tau,\tau]} ||u_N-w_N||_{H^{m}(\R^d_{x'})} = O(h^{N+1-m}).
    \]
    for integer $m$. And the coefficients of this expansion satisfy the conditions given in Theorem \ref{thm.main} (\eqref{eq.initial_conditions_abel} and either \eqref{eq.transport_sym_abel} or \eqref{eq.transport_anti_abel}).
    \par
    Similarly, $v_N$ has an expansion of the form
    \[
    v_N(h,V,W,x) = \varphi_V(x)\vec{e_l} + h\sum_{\pm m=0}^N e^{\frac{im}{h}\mscal{x,W}} \begin{bmatrix} \sum_{p=0}^{N-1} h^pB_{m,p}^{(1)}(x)\\  \sum_{p=0}^{N-1} h^pB_{m,p}^{(2)}(x)\\ ...\\ \sum_{p=0}^{N-1} h^p B_{m,p}^{(n)}(x)\end{bmatrix}
    \]
    and is unique up to order $N+1$; i.e. if $\tilde{w}_N$ satisfies \eqref{eq.initial_conditions_nonabel}, then
    \[
    \sup_{x_0 \in [-\tau,\tau]} ||v_N-\tilde{w}_N||_{H^m(\R^d_{x'})} = O(h^{N+1}).
    \]
    And the coefficients of $v_N$ satisfy the conditions given in Theorem \ref{thm.non_abelian} (\eqref{eq.initial_conditions_nonabel},\eqref{eq.transport_nonabel}).
\end{theorem}
\subsection{Proof of Theorems \ref{thm.main} and \ref{thm.non_abelian}}
If we have the approximate solutions given in \ref{thm.approx}, the last step remaining is to show that $u_N,v_N$ approximate exact solutions $u,v$ which satisfy Theorems \ref{thm.main} and \ref{thm.non_abelian}. To achieve this, we will use a result which applies to high-frequency oscillating solutions of a general class of hyperbolic systems. 
To state the theorem, we must introduce some semiclassical Sobolev spaces (see \cite{Zworski} for more details), which will correspond to the conditions \eqref{eq.exact_error_bounds_abel},\eqref{eq.approx_error_est}  on our solutions. 
\par
For $\Omega = [T_0,T] \times \R^d \subseteq \R^{1+d}$, $\rho > 0$, $\epsilon > 0$ and $m \in \mathbb{N}$, we have $w_h \in \mathcal{B}_{\rho,\epsilon}^m(\Omega)$ if for $h \in (0,\epsilon]$
\begin{gather*}
    w_h, \nabla w_h \in C^0([T_0,T],H^m(\R^d_{x'})) \text{ and }\\
    \sup_{x_0 \in [T_0,T]} ||(h\nabla_{x'})^k w_h(x_0,\cdot)||_{H^1(\R^d)} + \sup_{x_0 \in [T_0,T]} ||\partial_{x_0}(h\nabla_{x'})^k w_h(x_0,\cdot)||_{L^2(\R^d)} \leq \rho \\\text{ for all $0\leq k \leq m$.}
\end{gather*}
And $w_h \in \mathcal{A}_{\rho,\epsilon}^m(\Omega)$ if for $h \in (0,\epsilon]$
\begin{gather*}
    w_h,\nabla w_h \in C^0([T_0,T],W^{m,\infty}(\R^d_{x'})), \text{ and for all $1 \leq k \leq m$,}\\
    \sup_{x_0 \in [T_0,T]} ||w_h(x_0,\cdot)||_{W^{1,\infty}(\R^d)} \leq \rho, \sup_{x_0,\in [T_0,T]} ||(h\nabla_{x'})^k w_h(x_0,\cdot)||_{L^\infty(\R^n)} \leq \rho h, \text{ and}\\
    \sup_{x_0\in [T_0,T]} ||(h\nabla_{x'})^k \nabla w_h(x_0,\cdot)||_{L^\infty(\R^d)} \leq \rho.
\end{gather*}

\begin{theorem}[Gu\`es]\label{thm.Gues} Let integer $m > \frac{d}{2} + 1$, $M > m$, and $\Omega = [T_0,T] \times \R^d \subseteq \R^{1+d}$. For any $\rho > 0$, there exists $\epsilon_\rho > 0$ and $r > 0$ such that for any $w_h \in \mathcal{A}_{\rho,\epsilon_\rho}^{m+1}$ which satisfies 
\[
\square w_h = f(x,w_h,\nabla w_h) + h^M G_h(x)
\]
with $f(0,0,0) = 0$, $f$ smooth, and $G_h \in \mathcal{B}_{\rho,1}^m(\Omega)$, there exists a unique $u_h$ for $h \in (0,\epsilon_\rho]$ such that $u_h-w_h \in h^M \mathcal{B}_{r,\epsilon_\rho}^m(\Omega)$, and
\begin{gather}
    \square u_h = f(x,u_h,\nabla u_h)\\
    (u_h,\partial_{x_0}u_h)|_{x_0=T_0} = (w_h,\partial_{x_0}w_h)|_{x_0=T_0}.
\end{gather}
\end{theorem}
This is a particular application of Gu\`es's result in \cite{Gues}. A proof can be found in \cite{NSB}.
\par
Applying the theorem, we see that there exists a unique exact solution $u$ to \eqref{eq.system} satisfying the initial condition $u=u_{inc}$, and that $u-u_N = h^{N+1} \mathcal{E}_N(h,V,W,x)$, where $\mathcal{E}_N$ satisfies the estimates given in \eqref{eq.exact_error_bounds_abel}. This $u$ is the solution that meets the requirements of Theorem \ref{thm.main}.
\par
Similarly, this result tells us that there exists a unique exact solution $v$ to \eqref{eq.system} satisfying the initial condition $v= v_{inc}$, and that $v-v_N = h^{N+1} \tilde{\mathcal{E}}_N(h,V,W,x)$, where $\tilde{\mathcal{E}}_N$ satisfies \eqref{eq.exact_error_bounds_abel}. So $v$ is the solution that meets the requirements of Theorem \ref{thm.non_abelian}.
\subsection{Calculation of $(\square-Q_T)u_N$, $(\square-Q_T)v_N$}
    We begin with a list of calculations. First,
    \begin{equation*}
    \square \left(e^{\frac{im}{h} \mscal{x,W}} \begin{bmatrix}  A_m^{(1)}(x)\\  A_m^{(2)}(x)\\ ...\\  A_m^{(n)}(x) \end{bmatrix}\right) = 
    e^{\frac{im}{h} \langle x,W\rangle_M} \left(\frac{2im}{h} \T_W  \begin{bmatrix} A_m^{(1)}(x)\\  A_m^{(2)}(x)\\ ...\\ A_m^{(n)}(x)\end{bmatrix} + \square  \begin{bmatrix} A_m^{(1)}(x)\\  A_m^{(2)}(x)\\ ...\\  A_m^{(n)}(x)\end{bmatrix}\right)
\end{equation*}
Next, we expand our coefficients in $h$ as $ A_m^{(r)}= \sum_{p=0}^{N} h^pA_{m,p}^{(r)}$. Calculating the null forms for these expansions results in 
\begin{multline*}
Q_0\left(\sum_{\pm m = 0}^N A_m^{(r)} e^{\frac{im}{h} \mscal{x,W}},\sum_{\pm m = 0}^N A_m^{(r)} e^{\frac{im}{h} \mscal{x,W}}\right) = h^{-1}\sum_{\pm m = 0}^{2N} \sum_{\substack{\mu+\nu = m\\ |\mu|,|\nu| \leq N}} e^{\frac{im}{h} \mscal{x,W}}\Big(hQ_0(A_\mu^{(r)},A_\nu^{(r)}) \\
+R_{0,\mu,\nu}(A_\mu^{(r)},A_\nu^{(r)})\Big),
\end{multline*}
\begin{multline*}
    Q_0\left(\varphi_V,\sum_{\pm m = 0}^N A_m^{(r)} e^{\frac{im}{h} \mscal{x,W}}\right) = h^{-1} \sum_{\pm m = 0}^N e^{\frac{im}{h} \mscal{x,W}} \left(P_{0,m}(A_m^{(r)}) + h S_{0,m}(A_m^{(r)})\right),
\end{multline*}
\begin{multline*}
    Q_{ij}\left(\sum_{\pm m = 0}^N A_m^{(k)} e^{\frac{im}{h} \mscal{x,W}},\sum_{\pm m = 0}^N A_m^{(l)} e^{\frac{im}{h} \mscal{x,W}}\right) = \sum_{\pm m = 0}^{2N} \sum_{\substack{\mu+\nu = m\\ |\mu|,|\nu| \leq N}} h^{-1} e^{\frac{im}{h} \mscal{x,W}}\Big(hQ_{ij}(A_\mu^{(k)},A_\nu^{(l)}) \\
    + R_{ij,\mu,\nu}(A_\mu^{(k)},A_\nu^{(l)})\Big),
\end{multline*}
and
\begin{multline*}
     Q_{ij}\left( \varphi_V,\sum_{\pm m = 0}^N A_m^{(r)} e^{\frac{im}{h} \mscal{x,W}}\right) = h^{-1} \sum_{\pm m = 0}^N e^{\frac{im}{h} \mscal{x,W}} \left(P_{ij,m}(A_m^{(r)}) + hS_{ij,m}(A_m^{(r)})\right)
\end{multline*}
where
\begin{gather*}
    P_{0,m}(A_{m}^{(r)}) = im\varphi_V'(x) \mscal{V,W} A_m^{(r)}\\
    S_{0,m}(A_m^{(r)}) = \varphi_V'(x) \T_{V} A_m^{(r)}\\
    R_{0,\mu,\nu}(A_\mu^{(r)},A_\nu^{(r)}) =i\left(\mu A_\mu^{(r)}\T_W A_\nu^{(r)} + \nu A_\nu^{(r)} \T_WA_\mu^{(r)}\right)\\
    P_{ij,m}(A_m^{(r)}) = im\varphi_V'(x) (\theta_i \omega_j - \theta_j \omega_i) A_m^{(r)}\\
    S_{ij,m}(A_m^{(r)}) = \varphi_V'(x) (\theta_i \partial_j A_m^{(r)} - \theta_j \partial_i A_m^{(r)})\\
    R_{ij,\mu,\nu}(A_\mu^{(k)},A_\nu^{(l)}) = i \left(\mu A_\mu^{(k)}(\omega_i \partial_j A_\nu^{(l)} - \omega_j \partial_i A_\nu^{(l)}) + \nu A_\nu^{(l)}(\omega_j \partial_i A_\mu^{(k)} - \omega_i \partial_j A_\mu^{(k)})\right)
\end{gather*}
We also remind ourselves that  $Q_0(h^{a_k}\varphi_V,h^{a_l}\varphi_V) = Q_{jk}(h^{a_k}\varphi_V,h^{a_l}\varphi_V) = 0$, $Q_0(C,w) = Q_{ij}(C,w) = 0$, $Q_0(hu_1+u_2,w) = hQ_0(u_1,w) + Q_0(u_2,w)$ and $Q_{ij}(hu_1+u_2,w) = hQ_{ij}(u_1,w) + Q_{ij}(u_2,w)$.
\par
Applying the above to our approximate solutions, we obtain
\begin{multline*}
    Q_0(u_N^{(l)},u_N^{(r)}) = (1+\mathds{1}_{r=l})h^{2} \left(\sum_{\pm m = 0}^N e^{\frac{im}{h}\mscal{x,W}}(P_{0,m}(A_m^{(r)}) + hS_{0,m}(A_m^{(r)}))\right) + \\h^{3}\sum_{m=0}^{2N} \sum_{\substack{\mu + \nu = m\\ |\mu|,|\nu| \leq N}} e^{\frac{im}{h}\mscal{x,W}} \left(hQ_0(A_\mu^{(l)},A_\nu^{(r)})+R_{0,\mu,\nu}(A_\mu^{(l)},A_\nu^{(r)})\right)
\end{multline*}
and for $s \neq l$,
\begin{equation*}
    Q_0(u_N^{(s)},u_N^{(r)}) = h^{3}\sum_{m=0}^{2N} \sum_{\substack{\mu + \nu = m\\ |\mu|,|\nu| \leq N}} e^{\frac{im}{h}\mscal{x,W}} \left(hQ_0(A_\mu^{(s)},A_\nu^{(r)})+R_{0,\mu,\nu}(A_\mu^{(s)},A_\nu^{(r)})\right).
\end{equation*}
For $Q_{ij}$, if $r \neq l$,
\begin{multline*}
    Q_{ij}(u_N^{(l)}, u_N^{(r)}) =h^{2} \sum_{\pm m = 0}^N e^{\frac{im}{h} \mscal{x,W}} (P_{ij,m}(A_m^{(r)}) + hS_{ij,m}(A_m^{(r)})) \\
    + h^{3} \sum_{\pm m =0}^{2N} \sum_{\substack{\mu+\nu = m\\ |\mu|,|\nu| \leq N}} e^{\frac{im}{h}\mscal{x,W}} \left(hQ_{ij}(A_\mu^{(l)},A_\nu^{(r)})+R_{ij,\mu,\nu}(A_\mu^{(l)},A_\nu^{(r)})\right)
\end{multline*}
and if $s,r \neq l$ with $s \neq r$,
\begin{equation*}
    Q_{ij}(u_N^{(s)},u_N^{(r)}) = h^{3} \sum_{\pm m =0}^{2N} \sum_{\substack{\mu+\nu = m\\ |\mu|,|\nu| \leq N}} e^{\frac{im}{h}\mscal{x,W}} \left(hQ_{ij}(A_\mu^{(s)},A_\nu^{(r)})+R_{ij,\mu,\nu}(A_\mu^{(s)},A_\nu^{(r)})\right).
\end{equation*}
Lastly, due to anti-symmetry, $Q_{ij}(u_N^{(r)}, u_N^{(r)}) = 0$ for all $r$.
\par
From this, we can look at \eqref{eq.system}, and see that
\begin{multline}\label{eq.calc_result_abel}
    h^{2}\sum_{\pm m = 0}^N e^{\frac{im}{h}\mscal{x,W}} \left(\frac{2im}{h} \T_W A_m^{(r)} + \square A_m^{(r)}\right) = 
 \sum_{\pm m = 0}^N e^{\frac{im}{h} \mscal{x,W}} \\\sum_{\alpha=1}^n 2\Bigg(q_{r\alpha l}h^{2}(P_{0,m}(A_m^{(\alpha)}) + hS_{0,m}(A_m^{(\alpha)}))+ \sum_{\substack{i,j=1\\i > j}}^d  h^{2} p_{rijl\alpha}(P_{ij,m}(A_m^{(\alpha)}) + h S_{ij,m}(A_m^{(\alpha)})) \Bigg)\\
    +  \sum_{\pm m = 0}^{2N}e^{\frac{im}{h}\mscal{x,W}}\sum_{\substack{\mu+\nu = m\\|\mu|,|\nu| \leq N}} \sum_{\alpha,\beta = 1}^n\Bigg( h^{3} q_{r\alpha \beta} \Big(R_{0,\mu,\nu}(A_\mu^{(\alpha)}, A_\nu^{(\beta)})+hQ_0(A_\mu^{(\alpha)},A_\nu^{(\beta)})\Big)\\+
    \sum_{\substack{i,j=1\\i > j}}^d h^{3} p_{rij\alpha\beta}\Big(R_{ij,\mu,\nu}(A_\mu^{(\alpha)},A_\nu^{(\beta)})+hQ_{ij}(A_\mu^{(\alpha)},A_\nu^{(\beta)})\Big)\Bigg)
\end{multline}

\par
We similarly have
\begin{multline}\label{eq.calc_result_nonabel}
    \sum_{\pm m = 0} e^{\frac{im}{h} \mscal{x,W}}\left(2im \T_W B_m^{(r)} + h\square B_m^{(r)}\right) = 
    \sum_{\pm m = 0} e^{\frac{im}{h} \mscal{x,W}} \sum_{\alpha = 1}^n2\Big(q_{r\alpha l}(P_{0,m}(B_m^{(l)}) + hS_{0,m}(B_m^{(l)}))\\
    + \sum_{\substack{i,j = 1\\ i > j}}^d  p_{rijl\alpha} (P_{ij,m}(B_m^{(\alpha)}) + hS_{ij,m}(B_m^{(\alpha)}))\Big)\\
    + \sum_{\pm m = 0}^{2N} e^{\frac{im}{h} \mscal{x,W}} \sum_{\substack{\mu + \nu = m\\|\mu|,|\nu| \leq N}} \sum_{\alpha,\beta = 1}^n\Bigg(h q_{r\alpha\beta}(R_{0,\mu,\nu}(B_\mu^{(\alpha)},B_\nu^{(\beta)})+hQ_0(B_\mu^{(\alpha)},B_\nu^{(\beta)}))\\
    + \sum_{\substack{i,j=1\\i>j}}^d  hp_{rij\alpha\beta}(R_{ij,\mu,\nu}(B_\mu^{(\alpha)},B_\nu^{(\beta)}) + hQ_{ij}(B_\mu^{(\alpha)},B_\mu^{(\beta)})\Bigg)
\end{multline}
\par
\subsection{The Linear Regime}
If we group by $m$, \eqref{eq.calc_result_abel} becomes
\begin{equation}\label{eq.coeff_abel}
2im \T_W A_m^{(r)} = h\Big(-\square A_m^{(l)} + 2\sum_{\alpha=1}^nq_{r\alpha l} P_{0,m}(A^{(\alpha)}_m)+2\sum_{\substack{i,j=1\\i>j}}p_{rijl\alpha}P_{ij,m}(A_m^{(\alpha)})\Big)\Big) + O(h^2)
\end{equation}
Consider the expansion of $A_m^{(r)}$ in $h$ as $A_m^{(r)} = A_{m,0}^{(r)} + hA_{m,1}^{(r)}+...+h^NA_{m,N}^{(r)}$. Also, since $q_{r\alpha l},p_{rijl\alpha}$ are smooth, we can use the Taylor expansion $q_{r\alpha l}(x,u) = q_{r\alpha l}(x,C)+ h (x-C) \cdot (\nabla_uq_{r\alpha l})(x,C)+h^2 \tilde{q_{r\alpha l}}(h,V,W,x)$. Then if we make the initial assumption that \eqref{eq.initial_conditions_abel} defines $A_{m,0}^{(r)}$ for all $x$, we can group by the order of $h$ and solve for all $A_{m,p}^{(r)}$. 
\par
Since $A_{1,0}^{(r)} \equiv 0$ for $r \neq k$, if $k=l$ then this results in
\begin{equation}
    2i\T_W A_{1,1}^{(r)} = 2 q_{rll}(x,C) P_{0,1}(A_{1,0}^{(l)}).
\end{equation}
which is \ref{eq.transport_sym_abel}.
If instead $k \neq l$, we have
\begin{equation}
    2i\T_W A_{1,1}^{(r)} = 2q_{rkl}(x,C) P_{0,1}(A^{(k)}_{1,0})+2\sum_{\substack{i,j=1\\i>j}}p_{rijlk}(x,C) P_{ij,1}(A_{1,0}^{(k)}).
\end{equation}
which is \ref{eq.transport_anti_abel}.
\par
Solving for all $A_{m,p}^{(r)}$ is possible, since (for $m \neq 0$) the right side of the equation will be dependent on $A_{m,p-1}^{(l)}$ and the $O(h^2)$ term dependent on $A_{m',q}^{(r)}$ for $q < p-1$. Note that the initial condition for $x_0 \ll 0$ means these solutions are unique, and thus this will be the unique expansion of this form.
\begin{remark}
    For $m=0$, $2im\T_WA_0^{(r)} = 2P_{0,0}^{(l)}(A_0^{(l)}) = 0$, so we end up with a different equation
    \[
    \square A_m^{(l)} = O(h).
    \]
    But with our initial conditions, this also has a unique solution.
\end{remark}
\par
When we set $A_{m,p}$ to the coefficients solved for above, $(\square-Q_T)u_N$ consists of the terms which are the right-hand side of \eqref{eq.calc_result_abel} when $|m| > N$ or order $h > N$. However, since $A_{m,0}^{(r)} = 0$ for $|m| > 1$ and $Q,R$ are bilinear and show up as $O(h^3)$ in \eqref{eq.coeff_abel}, we see that $A_{m,p}^{(r)} = 0$ whenever $|m|>1$ and $p\leq |m|$.  So $(\square-Q_T)u_N$ consists of terms that are order $h > |N|$.
\par
Since all $A_{m,p}$ are smooth, this implies that we have \eqref{eq.approx_eval} with 
\[E_N = \sum_{N<|m|\leq 2N} e^{\frac{im}{h} \mscal{x,W}} G_{N,m}(h,V,W,x)\] where each $G_{N,m}$ is smooth in $x$ and polynomial in $h$. This implies that $E_N$ satisfies \eqref{eq.approx_error_est}.
\subsection{The weakly nonlinear regime}
Grouping \eqref{eq.calc_result_nonabel} by $m$ results in 
\[
    2im\T_W B_m^{(r)} =  2\sum_{\alpha = 1}^n \left(  q_{rl\alpha}P_{0,m}(B_m^{(\alpha)}) + \sum_{\substack{i,j = 1\\ i>j}} p_{rijl\alpha} P_{ij,m}(B_m^{(\alpha)})\right) + O(h).
\]
If we expand $q_{rkl}$,$p_{rijl\alpha}$ in $h$ as in the previous cases, with the initial assumptions \eqref{eq.initial_conditions_nonabel}, we obtain ODE systems which can be used to solve for each $B_{m,p}$. Of note, for $B_{1,0}$ we have
\[
2im\T_W B_{1,0}^{(r)} = 2\sum_{\alpha=1}^n \left(q_{rl\alpha}(x,C) P_{0,1}(B_{1,0}^{(l)}) + \sum_{\substack{i,j=1\\i>j}}  p_{rijl\alpha}(x,C)P_{ij,1}(B_{1,0}^{(\alpha)})\right)
\]
which is \eqref{eq.transport_nonabel}, as desired. We can follow the logic of the linear regime to obtain an approximate solution satisfying Theorem \ref{thm.approx}.

\end{document}